# A New Approach on the Seating Couples Problem

Daniel Kohen and Ivan Sadofschi

**Abstract**

A king invites $n$ couples to sit around a round table with $2n + 1$ seats. For each couple, the king decides a prescribed distance $d$ between 1 and $n$ which the two spouses have to be seated from each other (distance $d$ means that they are separated by exactly $d - 1$ chairs). We will show that there is a solution for every choice of the distances if and only if $2n + 1$ is a prime number using a theorem known as Combinatorial Nullstellensatz.

## 1 Introduction

We present a proof[1] of the king's problem different from the proof given in [1] by considering the same polynomial. This polynomial arises naturally when trying to use the Combinatorial Nullstellensatz. As the proof presented in [1], this one is nonconstructive.

The problem was originally asked in [4], where R. Bacher also conjectures that the following holds:

**Conjecture:** If every prescribed distance $d$ is invertible modulo $m = 2n + 1$, then the couples can sit around the table.

Unfortunately this method fails when considering this generalization, because it depends on the fact that $\mathbb{Z}/p\mathbb{Z}$ is a field.

---

[1]This article follows the ideas we presented in `http://grupofundamental.wordpress.com/2010/03/05/il-faut-exiger-de-chacun-ce-que-chacun-peut-donner/#comments`

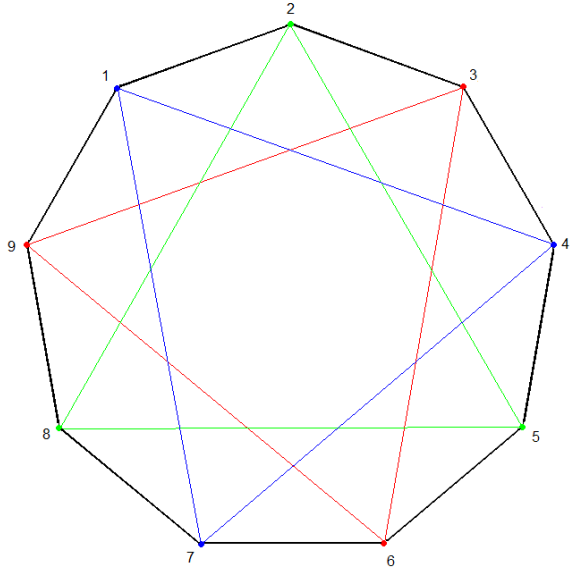

Figure 1: Orbits with $2n+1 = 9$ and $k = 3$

## 2 A counterexample in the composite case

First we shall prove that if $2n + 1$ is a composite number the king may choose the distances in a way that the couples cannot be seated. Let $k$ be a divisor of $2n + 1$ such that $1 < k < 2n + 1$.

Now let the distance to be $k$ for every couple. Number the seats from 1 to $2n + 1$, clockwise in a natural way. Let the orbit of $j$ be the seats with numbers congruent to $j$ modulo $k$.

Notice that since $k \mid 2n+1$, each orbit has the same number of seats, $\frac{2n+1}{k}$. Moreover, if two seats are at distance $k$ then they are in the same orbit. Since the number of seats in each orbit is odd, at least one seat per orbit must be empty. But there are $k > 1$ orbits, so the couples cannot seat around the table.

## 3 A solution in the prime case

Now suppose that $2n + 1$ is a prime number. We shall construct a polynomial in $\mathbb{F}_p[x_1, ..., x_n]$ having in mind the polynomial method. Next we show that this polynomial yields a solution to the king's problem using the Combinatorial Nullstellensatz.

In addition we shall make use of Dyson conjeture to calculate one of the coefficients of the polynomial. For the sake of completness we include the statements ot these two facts:

**Combinatorial Nullstellensatz:** Let $F$ be an arbitrary field, and let $f = f(x_1, \ldots, x_n)$ be a polynomial in $F[x_1, \ldots, x_n]$. Suppose that $\deg(f)$ is $\sum_{i=1}^{n} t_i$ where each $t_i$ is a nonnegative integer, and suppose the coefficient of $\prod_{i=1}^{n} x_i^{t_i}$ in $f$ is nonzero. Then, if $S_1, \ldots, S_n$ are subsets of $F$ with $|S_i| > t_i$, there are $s_1 \in S_1, s_2 \in S_2, \ldots, s_n \in S_n$ so that $f(s_1, s_2, \ldots, s_n) \neq 0$.

**Dyson Conjecture[2]:** The constant term in the Laurent Polynomial $\prod_{1 \leq i \neq j \leq n}(1 - \frac{x_i}{x_j})^{a_i}$ equals $\frac{(a_1 + a_2 + \ldots + a_n)!}{a_1! a_2! \ldots a_n!}$.

Now we are ready to proceed with the proof.

Given $n$ distances $d_1, \ldots, d_n$, a solution to the king's problem is a tuple $(x_1, \ldots, x_n)$ such that the $2n$ numbers $x_1, \ldots, x_n, x_1 + d_1, \ldots, x_n + d_n$ are pairwise distinct modulo $p = 2n + 1$.

Consider the following polynomial in $\mathbb{F}_p[x_1, ..., x_n]$:

$$f(x_1, \ldots, x_n) = \prod_{1 \leq i < j \leq n} (x_i - x_j)(x_i + d_i - x_j)(x_i - x_j - d_j)(x_i + d_i - x_j - d_j).$$

This polynomial takes values distinct from 0 if and only if it is evaluated at solutions to the kings problem, for it is the product of all pairwise differences of $x_1, ..., x_n, x_1 + d_1, \ldots, x_n + d_n$.

We have $\deg(f) = 4\binom{n}{2} = n(2n-2)$, so $m = x_1^{2n-2} \ldots x_n^{2n-2}$ is a monomial of degree $\deg(f)$.

To find the coefficient of $m$ in $f$, note that by taking only maximum degree terms, $f$ simplifies to

$$\prod_{1 \leq i < j \leq n}(x_i - x_j)^4 = x_1^{2n-2} \ldots x_n^{2n-2} \prod_{1 \leq i \neq j \leq n}(1 - \frac{x_i}{x_j})^2.$$

---
[2]Actually, a theorem proved independently by Wilson and Gunson in 1962. For a proof, see [3]

By Dyson's Conjecture, the constant term in $\prod_{1 \leq i \neq j \leq n} (1 - \frac{x_i}{x_j})^2$ equals $\frac{(2n)!}{2^n}$.

Now $(2n)! = (p-1)! \equiv -1 \mod p$ by Wilson's theorem.

In addition, by Fermat's Little Theorem: $p \mid 2^{p-1} - 1 = 2^{2n} - 1 = (2^n + 1)(2^n - 1)$ so $2^n \equiv \pm 1 \mod p$.

Therefore, we conclude that $[m]f = \frac{(2n)!}{2^n} \equiv \pm 1 \mod p$, so by the Combinatorial Nullstellensatz, there is a solution to the king's problem in this case, as we wanted to prove.